\documentclass[twocolumn,nofootinbib]{revtex4}

\makeatletter
\def\p@subsection{}
\makeatother

\usepackage{latexsym}
\usepackage{amssymb,amsmath}
\usepackage{amsthm}
\usepackage{nicefrac}
\usepackage{empheq}
\usepackage{multirow}
\usepackage{amssymb}
\usepackage{times}
\usepackage{amssymb}
\usepackage[np]{numprint}
\npdecimalsign{.}
\usepackage[utf8]{inputenc}
\usepackage{subfig} 
\usepackage{fancyhdr}
\usepackage{url}
\usepackage{tikz}
\usepackage{hyperref}
\usepackage{color}
\usepackage[textsize=tiny]{todonotes}
\usepackage{array}
\newcolumntype{H}{>{\setbox0=\hbox\bgroup}c<{\egroup}@{}}

\usepackage{graphicx}

\tikzstyle{node}=[fill=none, draw=black, text=black, rectangle, inner sep=1pt, minimum size=15pt]

\newcommand{\maxent}{q}
\newcommand{\distrFamily}{\mathcal{G}}
\newcommand{\Ex}[1]{\ensuremath{E\!\left(#1\right)}}
\newcommand{\mcProcess}{\vec{X}}
\newcommand{\mcState}{\vec{x}}
\newcommand{\changeVector}{v}

\newcommand{\meanVector}{\vec{\mu}}
\newcommand{\mcProcY}{\vec{Y}}
\newcommand{\mcProcZ}{\vec{Z}}
\newcommand{\mcStateY}{\vec{y}}
\newcommand{\mcStateZ}{\vec{z}}
\newcommand{\DNA}{\mathit{DNA}}
\newcommand{\bbbn}{\mathbb{N}}

\newcommand{\maxDevPercErr}{|| \epsilon ||_\infty^\%}

\newcommand{\dualFunctionThreshold}{\delta_\Psi}

\newcommand{\wsMCM}{\emph{weighted sum MCM}}
\newcommand{\jMCM}{\emph{joint MCM}}
\newcommand{\MM}{\emph{MM}}

\theoremstyle{definition}
\newtheorem{exmp}{Example}
 
\begin{document}

\title{\bf Reconstruction of Multimodal Distributions for Hybrid Moment-based Chemical Kinetics\\
Supporting Information}
\author{Alexander Andreychenko$^1$}
\thanks{Corresponding author. Email: alexander.andreychenko@uni-saarland.de}
\author{Linar Mikeev$^1$}
\author{Verena Wolf$^1$ \vspace*{3mm}}
\affiliation{$^1$ Computer Science Department, Saarland University, 
Saarbrucken, Germany 66123}



\begin{abstract}
\noindent The stochastic dynamics of biochemical reaction networks can be accurately 
described by discrete-state Markov processes where each chemical reaction
corresponds to a state transition of the process.
Due to the largeness problem of the state space, analysis techniques based on
an exploration of the state space 
are often not feasible and the integration of the moments of the underlying
probability distribution has   become a
very popular alternative. 
In this paper the focus is on a comparison of   reconstructed
 distributions from their moments obtained by two different moment-based analysis methods, 
the method of moments (MM) and the method of conditional moments (MCM). 
We use the maximum entropy principle to  derive a distribution 
that fits best to a given sequence of (conditional) moments.
For the  two gene regulatory networks that we consider we find that
  the MCM approach is more suitable to describe multimodal
  distributions and that the reconstruction of marginal distributions
   is more accurate if conditional distributions
are considered.

\vspace*{2ex}\noindent\textit{\bf Keywords}: Chemical Master Equation, Moment Closure, Method of Conditional Moments, Maximum Entropy.
\\[3pt]
\\[3pt]
\end{abstract}

\maketitle

 \newpage
 \clearpage
 \begin{appendix}
\section{Supporting Information}
In the following  sections we 
first 
describe in detail how the moment equations are obtained 
(Section~\ref{appenix:method_of_moments} and~\ref{appendix:mcm}) 
and how we approximate the support of the distribution (Section~\ref{appendix:SupportApprx}). 
In Section~\ref{appendix:maxEntProblem2D}
we then discuss the differences that arise during  the reconstruction of distributions 
with two instead of only one dimension,
and in Section~\ref{appendix:reconstructionMethods} 
we discuss the details of the reconstruction for the distributions of  
 the case studies introduced before.
More numerical results for the  two case studies are then provided
in Section~\ref{appendix:case_studies_more_details}.

\subsection{Method of Moments}
\label{appenix:method_of_moments}
For the time derivative of the expectation of a function 
$f:\mathbb{N}^n_0\to \mathbb{R}^n$ applied to
the vector of species, we directly get from 
Eq. (1)
	\begin{equation}\label{eq:expgen}
	\begin{array}{l}
	\frac{d}{dt}  \Ex{f(\mcProcess(t))} 
	 = 
	\sum\limits_{x} f(x) \frac{d}{dt} p(\mcState,t) \\
	\hspace{2ex} =
	 \sum\limits_{j=1}^m \Ex{\alpha_j(\mcProcess(t)) 
	(f(\mcProcess(t)\! +\!   \changeVector_j )\! -\! f(\mcProcess(t)))}.
	\end{array}
	\end{equation}
For $f(\mcState)=\mcState$ this yields a system of equations for the population means
\begin{equation}\label{eq:mean1}
\begin{array}{lcl}
  \frac{d}{dt}\Ex{\mcProcess(t)} &=& \sum\limits_{j=1}^m \changeVector_j \Ex{\alpha_j(\mcProcess(t))}.
\end{array}
\end{equation}
Note that the system of ODEs in Eq.~\eqref{eq:mean1} is only closed
if at most monomolecular reactions ($\sum_{i=1}^n\ell_{j,i}\le 1$) are involved.
Otherwise $\Ex{\alpha_j(\mcProcess(t))}$ involves moments of 
the
second order.
However, in this case we can approximate the unknown second order moments,
say $\Ex{X_i(t)\cdot X_{i'}(t)}$ if the reaction is of the form $S_i+S_{i'}\to\ldots$, $i\neq i'$,
either by assuming that the covariance is zero, which gives
$\Ex{X_i(t)\cdot X_{i'}(t)}=\Ex{X_i(t)}\cdot \Ex{X_{i'}(t)}$ or by
extending the system in~\eqref{eq:mean1}  with additional
equations for the second moments.
The general strategy is to replace $\alpha_j(\mcProcess(t))$ by a Taylor series about
the mean $\Ex{\mcProcess(t)}$.
Let us write $\mu_i(t)$ for $\Ex{X_i(t)}$ 
and $\meanVector(t)$ for the vector with entries $\mu_i(t)$, $1\le i\le n$. Then
	\begin{equation}\label{eq:TaylorExp}
	\begin{aligned}
	 & \hspace{-1.7ex}\Ex{\alpha_j(\mcProcess)} = \alpha_j(\meanVector) + \frac{1}{1!}  
	\textstyle\sum\limits_{i=1}^n\Ex{X_i-\mu_i}\frac{\partial}{\partial x_i}\alpha_j(\meanVector)\\[1ex]
	   & \hspace{-1.7ex}+  \frac{1}{2!} \! \textstyle\sum\limits_{i=1}^n \sum\limits_{k=1}^n \! 
	 \Ex{(X_i\!-\!\mu_i)(X_k\!-\!\mu_k)} 
	 \frac{\partial^2}{\partial x_i \partial x_k}\alpha_j(\meanVector)\\[1ex]
	 &\hspace{-1.7ex}+ \!\textstyle\ldots
	\end{aligned}
	\end{equation}
 where we omitted $t$ in the equation to improve readability.
 Note that  
 $\Ex{X_i(t)-\mu_i}=0$ and since we restrict 
 to reactions that are at most bimolecular
 with mass action kinetics, 
 all terms of order three or more disappear.
 The derivation of moments
 for general kinetics is 
 presented
 in~\cite{StumpfJournal}.
 
 By letting $C_{ik}$ be the covariance
 $\Ex{(X_i(t)-\mu_i)(X_k(t)-\mu_k)}$, we get
 \begin{equation}\label{eq:meanrate}
 \begin{array}{rcl}
	\Ex{\alpha_j(\mcProcess)} \!=\! \alpha_j(\meanVector) 
	+ \frac{1}{2} \sum\limits_{i=1}^n \sum\limits_{k=1}^n 
	C_{ik}\frac{\partial^2}{\partial x_i \partial x_k}\alpha_j(\meanVector).
 \end{array}
 \end{equation}
 Next, we derive an equation for the covariances by first exploiting the relationship
	\begin{equation}\label{eq:rel2nd}
	\begin{array}{l}
 \frac{d}{dt}  C_{ik} =  \frac{d}{dt} \Ex{X_iX_k}- \frac{d}{dt}(\mu_i\mu_k)\\[1ex] 
   =  \frac{d}{dt} \Ex{X_iX_k}- \left(\frac{d}{dt}\mu_i\right)\mu_k
  -\mu_i \left( \frac{d}{dt}\mu_k \right)
	\end{array}
	\end{equation}
 and if we couple this equation with the equations for the means, the only unknown term
 that remains is the derivative $\frac{d}{dt} \Ex{X_iX_k}$ of the second moment.
 For this we can use the same strategy as before, i.e., from Eq.~\eqref{eq:expgen}
 we get
\begin{equation}\label{eq:secondm}
	\begin{array}{l}
  \frac{d}{dt}  \Ex{X_i X_k} \! =\!
  \sum\limits_{j=1}^m \bigg( v_{j,i}v_{j,k}  \Ex{\alpha_j(\mcProcess)} \\[1ex]
  \hspace{2ex}+ v_{j,k}\Ex{\alpha_j(\mcProcess)X_i}
 + v_{j,i}\Ex{\alpha_j(\mcProcess)X_k}\bigg),
	\end{array}
\end{equation}
 where $v_{j,i}$ and $v_{j,k}$ are the corresponding entries of the vector $v_j$.
 Clearly, we can use Eq.~\eqref{eq:meanrate} for the term 
 $E(\alpha_j(\mcProcess))$, while
 the terms $E(\alpha_j(\vec X)X_i)$ and $E(\alpha_j(\mcProcess)X_k)$ have to be
 replaced by the corresponding Taylor series about the mean.
 Let $f_j(\mcState):=\alpha_j(\mcState)x_i$.
 Similar to Eq.~\eqref{eq:meanrate}, we get
 that $E(\alpha_j(\vec X)X_i)$ equals
   \begin{equation}\label{eq:meanrate2}
	\begin{aligned}[c]
& \textstyle\hspace{-1.5ex}\alpha_j(\meanVector) \mu_i+ \frac{1}{1!}  \sum\limits_{i=1}^n\Ex{X_i-\mu_i}\frac{\partial}{\partial x_i}f_j(\meanVector)\\[1ex]
&\textstyle \hspace{-1.5ex}+ \frac{1}{2!} \sum\limits_{i=1}^n \sum\limits_{k=1}^n\Ex{(X_i-\mu_i)(X_k-\mu_k)} \frac{\partial^2}{\partial x_i \partial x_k}f_j(\meanVector)\\[1ex]
& \textstyle\hspace{-1.5ex}+\ldots
	\end{aligned}
  \end{equation}
 
  Here, it is important to note that moments of order three come into play, since derivatives
 of order three of $f_j(\mcState)=\alpha_j(\mcState)x_i$ may be nonzero.
 It is possible to take these terms into  account
 by deriving additional equations for moments of order three and higher.
 Obviously, these equations will then include moments of even higher order such that
 theoretically we end up with an infinite system of equations. However, a popular
 strategy is to close the equations by assuming that all moments of order  $>M$
  that are centered around the mean are equal to zero. E.g. if we choose $M=2$, then
 we can simply use the approximation
 \begin{equation*}
 	\begin{aligned}[c]
 &\textstyle\Ex{\alpha_j(\mcProcess)X_i} \approx 
  \alpha_j(\meanVector) \mu_i \\
 & \textstyle+ \frac{1}{2!} \sum\limits_{i=1}^n \sum\limits_{k=1}^n\Ex{(X_i-\mu_i)(X_k-\mu_k)} \frac{\partial^2}{\partial x_i \partial x_k}f_j(\meanVector).
 	\end{aligned}
 	\end{equation*}
 
 Other methods can be used to close the system of equations,
 including 
 derivative matching and zero cumulants closure~\cite{hespanha_moment_2008},
 as well as those that make
 assumptions about the distribution of 
 the moments~\cite{gomez-uribe_mass_2007,milner_moment_2011,pubsdoc:moment-closures-non-linear-rates}.
 
 The given approximation is then inserted into Eq.~\eqref{eq:secondm}
 and the result is used to replace  the term 
$\frac{d}{dt} \Ex{X_iX_k}$ in Eq.~\eqref{eq:rel2nd}.
 Finally, we can integrate the time evolution of the means and that of the covariances
 and variances.
 
 \stepcounter{exmp}
 \begin{exmp}\label{ex:geneExpression2}
 	We apply the moment closure technique described above to the
 	gene expression system from Example 1. 
 	When we consider only the moments up to second order,
 	the corresponding  equations for the average number of molecules are, for instance,
	given by
	\begin{equation*}
	\begin{aligned}[c]
		& \frac{d}{dt} \mu_{D_{\mbox{\scriptsize off}}} 
			 = \tau_{\mbox{\scriptsize off}} \mu_{D_{\mbox{\scriptsize on}}} 
			- \Ex{\tau_{\mbox{\scriptsize on}}^p X_{D_{\mbox{\scriptsize off}}} X_P}\\
		& \frac{d}{dt}\mu_{D_{\mbox{\scriptsize on}}} 
			 = \tau_{\mbox{\scriptsize on}} \mu_{D_{\mbox{\scriptsize off}}}
			+ \Ex{\tau_{\mbox{\scriptsize on}}^p X_{D_{\mbox{\scriptsize off}}} X_P} \\
		&	\frac{d}{dt}\mu_R = k_r \mu_{D_{\mbox{\scriptsize on}}} - \gamma_r \mu_R \\
		&	\frac{d}{dt}\mu_P = k_p \mu_R - \gamma_p \mu_P,
		\end{aligned}
	\end{equation*}	
	where $\mu_{D_{\mbox{\scriptsize off}}}, \mu_{D_{\mbox{\scriptsize on}}}$ 
	are the expected numbers of $D_{\mbox{\scriptsize off}}$ and $D_{\mbox{\scriptsize on}}$, 
	respectively, 
	and $\mu_R,\mu_P$ are the expected numbers of mRNA and proteins.
		
	Next we compare the obtained moments with those computed 
	via a direct numerical integration of the CME (Table~\ref{tab:gexpr_mc_error}). 
	We consider the following three cases.
	The moment closure approximation is
	carried out using all moments up to order 4, 6, and 8.
	For each case we list the number of moment equations, the running time, and
	the relative errors in the first four moments (columns 4-7).
	The relative error 
	for the moments of order $l$
	for species $i$
	 is given by
	$ 
		\epsilon^r_l = \max_{1 \leq i \leq n}
			\nicefrac
						{\vert \mu^{(i)}_l - \check\mu^{(i)}_l \vert}
						{\check\mu^{(i)}_l}$,
	where $\mu^{(i)}_l$ and $\check\mu^{(i)}_l$ are the values of the moments
	computed using the moment closure and obtained with a
	direct integration of the CME.
	
		Please note that in the reconstruction procedure we do 
		not use the moment of the highest order.
		For example, if we approximate moments up to order $6$,
		then the highest order that is taken into account during 
		the reconstruction is $5$ (which corresponds to the case $M=5$,
		cf. Sect.~\ref{appendix:case_studies_more_details})
		because of the high sensitivity of the numerical procedure
		even to the small absolute error in the moment approximation.
\begin{table}[b]
\caption{Moment closure approximation results for the gene expression system \label{tab:gexpr_mc_error}}
\centering
\begin{tabular}{cccccccc}
\hline
\shortstack{moment\\ closure\\ order} & $\#$ equ. &  \shortstack{time\\(sec)}& \shortstack{error\\ ord. 1\\ moments} & \shortstack{error\\ ord. 2\\ moments} & \shortstack{error\\ ord. 3\\ moments} & \shortstack{error\\ ord. 4\\ moments} \\
\hline
4 & 70 & 1 & \numprint{8e-6} & \numprint{8.3e-5} & \numprint{9.6e-5} &
\numprint{8.24e-4} \\
\hline
6 & 209 & 25 & \numprint{2e-6} & \numprint{2e-6} & \numprint{1e-5} & \numprint{3.6e-5} \\
\hline
8 & 494 & 3726 & \numprint{1e-6} &\numprint{2e-6} & \numprint{2e-6} &\numprint{4e-6}\\
\hline
\end{tabular}
\end{table}
\end{exmp}


\subsection{Method of Conditional Moments}
\label{appendix:mcm}
We first decompose the chemical populations described by $\mcProcess(t)$ into 
small and large  populations. Here we assume that this decomposition is static.
However, it is obvious that during the integration over time, we can (after reconstructing the joint distribution) choose a different 
decomposition for the remaining time.
From what size on a 
population should be considered as small is typically dependent on the amount of
main memory that is available and on the maximum order of the moments that we
consider for the large populations.
Note that considering   conditional moments 
yields a smaller number of equations if the   order of the considered moments is
high. The reason is that the number of equations for  representing the 
dynamics of the small populations does not increase as the order of considered 
conditional moments increases. Also, for many systems the decomposition is
obvious, as the small populations are exactly those that have a maximal size of, say, less than
10 (because they represent binding sites) and the large populations count 
protein 
numbers which may become rather large. 

Formally, we write the random vector 
$\mcProcess(t)$ at time $t$ as 
$\mcProcess(t)=(\mcProcY(t),\mcProcZ(t))$,
where $\mcProcY(t)$ corresponds to the small, and $\mcProcZ(t)$ to the large populations. 
Similarly, we write $\mcState = (\mcStateY,\mcStateZ)$ for the states of the process
and $\changeVector_j=(\hat{\changeVector}_j,\tilde \changeVector_j)$
for the change vectors, $j \in \{1,\ldots,m\}$.
Again, 
the first component refers to the small and the second component to the large populations.
Now, Eq. 1
becomes 

\begin{equation}
\label{eq:master2}
\hspace{-2ex}\begin{array}{l@{\,}c@{\,}r}
\frac{d p(\mcStateY,\mcStateZ)}{dt} 
&=& \sum\limits_{j=1}^m ( \alpha_j(\mcStateY\!-\!\hat \changeVector_j,
\mcStateZ \!-\! \tilde \changeVector_j)  p(\mcStateY\!-\!\hat 
\changeVector_j,\mcStateZ\!-\!\tilde \changeVector_j) \\[1ex]
 &&-  \alpha_j(\mcStateY,\mcStateZ)p(\mcStateY,\mcStateZ))
 \end{array}
\end{equation}
	where we omitted the time parameter $t$ to improve readability.
	Next, we sum over all possible $\mcStateZ$ to get the time evolution of the 
	marginal distribution $\hat p(\mcStateY)= \sum_{\mcStateZ} p(\mcStateY,\mcStateZ)$ 
	of the small populations.
	
	\begin{equation}\label{eq:smallcme}
	\begin{array}{l}
 \frac{d}{dt}\hat p(\mcStateY) =\\[1ex]
 \sum\limits_{\mcStateZ} \sum\limits_{j=1}^m \alpha_j(\mcStateY-\hat \changeVector_j,\mcStateZ-\tilde \changeVector_j)p(\mcStateY-\hat \changeVector_j,\mcStateZ-\tilde \changeVector_j) \\[2ex]
	\hspace{4ex} - \sum\limits_{\mcStateZ} \sum\limits_{j=1}^m \alpha_j(\mcStateY,\mcStateZ)p(\mcStateY,\mcStateZ)=\\[2ex]
  \sum\limits_{j=1}^m  \hat p(\mcStateY-\hat \changeVector_j) 
  E[\alpha_j(\mcStateY-\hat \changeVector_j,\mcProcZ)\mid Y=\mcStateY-\hat \changeVector_j]\\[2ex]
	\hspace{4ex}- 
	\sum\limits_{j=1}^m  \hat p(\mcStateY) 
	E[\alpha_j(\mcStateY,\mcProcZ)\mid \mcProcY=\mcStateY]
	\end{array}
	\end{equation}
	
	Note that in this small master equation that describes the change of the mode
	probabilities over time, the sum runs only over those reactions that modify 
	$\mcStateY$,
	since for all other reactions the terms cancel out. Moreover, on the right side
	we have only mode probabilities of neighboring modes and conditional
	expectations of the continuous part of the reaction rate. 
	For the latter, we can   use a Taylor expansion about the conditional
	population means. Similar to Eq.~\eqref{eq:TaylorExp}, this yields an 
	equation that involves the conditional means and centered conditional
	moments of second order (variances and covariances).
	Thus, in order to close the system of equations, we need to derive equations
	for the time evolution of the  conditional means and centered conditional
	moments of higher order. 
	Since the mode probability $p(\mcStateY)$ may 
	become zero, we first derive an equation for the evolution of the partial means
	(conditional means multiplied by the probability of the condition)
		\begin{equation*}
		\begin{array}{l}
			\frac{d}{dt}  \left(E[\mcProcZ \mid \mcStateY]\ p(\mcStateY) \right)
			= \sum\limits_{\mcStateZ} \mcStateZ \frac{d}{dt}  p(\mcStateY,\mcStateZ)\\[1ex]
		\hspace{2ex}	 = 
		\sum\limits_{j=1}^m   
			 E[(\mcProcZ + \tilde\changeVector_j) 
			 	\alpha_j(\mcStateY-\hat \changeVector_j, \mcProcZ)
			 	\mid \mcStateY-\hat \changeVector_j]\ 
			 	p(\mcStateY-\tilde \changeVector_j)\\[1ex]
			\hspace{2ex} - 
			\sum\limits_{j=1}^m   
			E[\mcProcZ \alpha_j(\mcStateY,\mcProcZ)	\mid \mcStateY]\ p(\mcStateY)
			,
		\end{array}
		\end{equation*}
	where in the  second line we   applied Eq.~\eqref{eq:master2} and simplified
	the result.
	The conditional expectations 
	$E[(\mcProcZ + \tilde\changeVector_j) 
				 	\alpha_j(\mcStateY-\hat \changeVector_j, \mcProcZ)
				 	\mid \mcStateY-\hat \changeVector_j]$
		and 
		$E[\mcProcZ \alpha_j(\mcStateY,\mcProcZ)	\mid \mcStateY]$
	are then replaced by their Taylor expansion about the conditional 
	means such that the equation involves only
	conditional means and  higher centered conditional moments~\cite{MCM_Hasenauer_Wolf}.
	For higher centered conditional moments, similar equations can be derived. If all
	centered conditional moments of order higher than $k$ are assumed to be zero,
	the result is a (closed) system of differential algebraic equations (algebraic equations are
	obtained whenever a mode probability 
	$p(\mcStateY)$ is equal to zero).
	However, it is possible to transform the system of differential algebraic equations 
	into a system of (ordinary) differential   equations
	after truncating modes with insignificant probabilities.
	Then we can get an accurate 
	approximation of the solution 
	after applying standard numerical integration methods.
	We construct the ODE system 
	using the tool 
	\setcounter{footnote}{0}
	SHAVE\footnote{L. Mikeev, \href{http://almacompute.mmci.uni-saarland.de/shave/}{http://almacompute.mmci.uni-saarland.de/shave/} }
	which implements 
	the truncation based approach
	and solve it using MATLAB's	\texttt{ode45} solver
	with the default error tolerance settings.
	  \begin{table}[b]
\caption{Conditional moment closure approximation results for the gene expression system \label{tab:gexpr_mcm_error}}
\centering
\resizebox{\columnwidth}{!}
{
\begin{tabular}{cccccccc}
\hline
\shortstack{cond.\\ moment\\ closure\\ order} & $\#$ eq. & \shortstack{time\\(sec)} & \shortstack{error\\ cond.\\ probs.} & \shortstack{error\\ ord. 1\\ cond.\\ moments} & \shortstack{error\\ ord. 2\\ cond.\\ moments} & \shortstack{error\\ ord. 3\\ cond.\\ moments} & \shortstack{error\\ ord. 4\\ cond.\\ moments} \\
\hline
4 & 30 & 1 & \numprint{7e-6} & \numprint{1e-5} & \numprint{2.86e-4} & \numprint{1.12e-3} & \numprint{6.98e-3} \\
\hline
6 & 56 & 2 & \numprint{6e-6} & \numprint{3.6e-5} & \numprint{5.9e-5} & \numprint{6.8e-5} & \numprint{2.18e-4} \\
\hline
8 & 90 & 9 & \numprint{2e-6} & \numprint{4.2e-5} & \numprint{6.2e-5} & \numprint{7.7e-5} & \numprint{9.1e-5} \\
\hline
\end{tabular}
}
\end{table}

	\begin{exmp}\label{ex:geneExpression_MCM}
		We apply the method of conditional moments to the
		gene expression system from Example~1. 
		The modes of the system are then given by the state of the DNA. 
		The equations for the mode probabilities ($p_{\mbox{\scriptsize off}}$, $p_{\mbox{\scriptsize on}}$) and the expected number of mRNA ($\mu_{R,\mbox{\scriptsize off}}$, $\mu_{R,\mbox{\scriptsize on}}$) and proteins ($\mu_{P,\mbox{\scriptsize off}}$, $\mu_{P,\mbox{\scriptsize on}}$) are as follows:
	\begin{equation*}
	\begin{aligned}[c]
&\frac{d}{dt} p_{\mbox{\scriptsize off}}  = \tau_{\mbox{\scriptsize on}} \: p_{\mbox{\scriptsize on}} - (\tau_{\mbox{\scriptsize off}} + \tau^p_{\mbox{\scriptsize on}} \: \mu_{P,\mbox{\scriptsize off}}) p_{\mbox{\scriptsize off}} \\
&\frac{d}{dt} \left( \mu_{R,\mbox{\scriptsize off}} \: p_{\mbox{\scriptsize off}} \right) = - \gamma_r \mu_{R,\mbox{\scriptsize off}} \: p_{\mbox{\scriptsize off}} \\
&\frac{d}{dt} \left( \mu_{P,\mbox{\scriptsize off}} \: p_{\mbox{\scriptsize off}} \right)  = ( k_p \mu_{R,\mbox{\scriptsize off}} - \gamma_p \mu_{P,\mbox{\scriptsize off}} ) p_{\mbox{\scriptsize off}}\\
& \frac{d}{dt} p_{\mbox{\scriptsize on}}  = ( \tau_{\mbox{\scriptsize off}} + \tau^p_{\mbox{\scriptsize on}} \mu_{P,\mbox{\scriptsize off}}) p_{\mbox{\scriptsize off}} - \tau_{\mbox{\scriptsize on}} \: p_{\mbox{\scriptsize on}} \\
& \frac{d}{dt} \left( \mu_{R,\mbox{\scriptsize on}} \: p_{\mbox{\scriptsize on}} \right) = ( k_r - \gamma_r \mu_{R,\mbox{\scriptsize on}} ) p_{\mbox{\scriptsize on}} \\
& \frac{d}{dt} \left( \mu_{P,\mbox{\scriptsize on}} \: p_{\mbox{\scriptsize on}} \right)  = ( k_p \mu_{R,\mbox{\scriptsize on}} - \gamma_p \mu_{P,\mbox{\scriptsize on}} ) p_{\mbox{\scriptsize on}} \\
	\end{aligned}
	\end{equation*}
	We computed the conditional moments and conditional probabilities 
	of the running example 
	(cf.\ Ex.~1 and Ex.~2) over time by
	considering moments up to the order of 4, 6, and 8. 
	For these three cases 
	the number of equations, when compared to the method of moments (MM), are as follows:
	\begin{center}
	{
		\begin{tabular}{lccc}
		\hline
		\ moment order $M$ & $\quad 4\quad $ &$\quad  6 \quad$ &$\quad  8\quad $ \\ \hline
		\ \#  equations for MM & $69$ & $209$ & $494$ 
		\\ \hline
		\ \#  equations for MCM\hspace{2ex} & $30$ & $56$ & $90$ 
		\\ \hline
		\end{tabular}
	}
	\end{center}
	The relative errors
	$\epsilon^r_l$ of the results of the method of conditional moments (MCM)  
	are given in Table~\ref{tab:gexpr_mcm_error},
	where we again compared to the 
	results obtained via a direct numerical solution. 


Our experiments show that the MCM performs much faster (due to the smaller number of equations) 
and still yields accurate approximation of the moments. 
For the chosen set of parameters
the MCM tends to provide a better approximation for higher moments,
whereas the MM approach is more accurate for lower moments when the same number of moments is considered.
For example, in the case of 6 moments
the maximum relative error for the first moments 
computed by the MM approach is 
$\numprint{2e-6}$,
compared to 
$\numprint{3.2e-5}$
when computed using the MCM.
At the same time, the maximum relative errors of the sixth moments 
are $\numprint{6.5e-4}$ and $\numprint{2e-4}$ for the MM and the MCM respectively. 
Note that the (unconditional) moments for the MCM are computed 
via multiplication of the conditional moments with the mode probabilities
and sum over all possible conditions.
We only consider non-central 
moments because the central moments introduce additional difficulties in the reconstruction
framework.

We also consider another set of parameters for the gene expression kinetics.
The rate constants are chosen 
$(\tau_{\mbox{\scriptsize on}},\tau_{\mbox{\scriptsize off}},k_r, k_p, 
\gamma_r, \gamma_p, \tau_{\mbox{\scriptsize \it on}}^p) = (0.05, 0.05, 10, 1, 4, 1, 0.015)$
as in~\cite{MCM_Hasenauer_Wolf}.
For the initial states we simply use
$x_{0,1}=(1,0,4,10)$ and $x_{0,2}=(1,0,4,10)$
with probabilities 
$P(x_{0,1}) = 0.7$ and $P(x_{0,1}) = 0.3$.
The comparison of the moment values at time instant
$t=10$
reveals that the MCM provides a much better approximation
both for high and low order moments
as opposed to the first parameter set.
For instance, in the case of 6 moments the
maximum relative error for the first moments
computed by the MM approach is $0.14$ whereas
in the MCM approach the error is 
$\numprint{7.5e-5}$.
The maximum relative error of the sixth moments
for the MM approach is $0.28$ compared
to $0.02$ using the MCM.
\end{exmp}
	
\subsection{Approximation of the Support}
\label{appendix:SupportApprx}
	During the iteration we  approximate the moments using 
	Eq. (6), where we do not sum
	over all states $x \in \bbbn_0$ but
	consider a subset $D = \{x_L,\ldots,x_R\}\subset \bbbn_0$. 
	Note that we have to find appropriate values for $x_L$ and $x_R$,
	since
	the iteration might fail to converge
	if the chosen value of $x_R$ is very large (and if $x_L=0$)
	as 
	the conditional number of the matrix 
	$ \left( H + \gamma^{(\ell)} \cdot \mathrm{diag}(H) \right) $
	is very large in this case.
		Thus, 
		we make use of the results in~\cite{Tari_asimplified} to find a region
		that contains the main part of the probability mass.
		We consider the roots of the   function  
		\begin{equation}
		\label{eq:delta0}
			\Delta^0(w) = 
					\left|
						\begin{array}{cccc}
							\mu_0 & \mu_1 & \cdots & \mu_k \\
							\vdots & &  & \vdots \\
							\mu_{k-1} & \mu_k & \cdots & \mu_{2k-1} \\
							1 & w & \cdots & w^k \\
						\end{array}
					\right|,
		\end{equation}
		where $k=\lfloor \frac{M}{2} \rfloor$, 
		and $M$ is even.
		Let $W = \{ w_1, \ldots, w_k \}$ be the set of the 
		solutions of   $\Delta^0(w) = 0$,
		where $w_1 < \ldots < w_k$ 
		are real and simple roots.
		The set   $D^{(0)} = \{ x_L^{(0)}, \ldots,  x_R^{(0)} \}$ with
	  $x_L^{(0)} = \lfloor w_1 \rfloor$ and $x_R^{(0)} = \lceil w_k \rceil$ 
	  is used as an initial guess
		for the approximated support when we start the optimization procedure. 
		The final results $ \lambda^*$ and $D^{*}$ of the iteration yields the distribution 
		\begin{center}
		$\tilde{q}(x) = \exp \left( -1 -\sum\nolimits_{k=0}^{M} \lambda^*_k x^k \right)$,
		\end{center} which 
		is an approximation 
		of the marginal distribution 
		$p_\centerdot(x,t) = P \left( X_\centerdot(t) = x \right)  $, i.e.\
		$p_\centerdot(x,t) \approx \tilde{q}(x) \mbox{ if } x \in D^{*}
		\mbox{ and } 
		p_\centerdot(x,t) \approx  0 \mbox { if } x \notin D^{*}.$
		\newline
	We can also
	account for the case of an odd number of moments.
	In addition to the function $\Delta^0(w)$ defined in Eq.~\eqref{eq:delta0},
	we also consider the function $\Delta^1(\eta)$
	\begin{equation*}
			\Delta^1(\eta) = 
				\left|
					\begin{array}{ccc}
							\mu_1 - \eta_1 \mu_0 & \cdots & \mu_{z} - \eta_1 \mu_{z-1} \\
							\vdots & \vdots  & \vdots \\
							\mu_{z-1} - \eta_1 \mu_{z-2} & \cdots & \mu_{2z-2} \! -\! \eta_1 \mu_{2z-3} \\
							1 & \cdots & \eta^{z-1} \\
					\end{array}
				\right|,
	\end{equation*}
	where $z=\lfloor \frac{M}{2} \rfloor + 1$ and $w_1$ is 
	the smallest root of the equation $\Delta^0(w)=0$.
	Again, let $W = \{w_1, \ldots, w_k \}$ be the set of the solutions
	of $\Delta^0(w)=0$ and $H= \{ \eta_1, \ldots, \eta_z \}$ be
	the set of solutions of $\Delta^1(\eta) = 0$,
	where all the elements of $W$ and $H$ are real and simple.
	The first approximation for the truncated support of the distribution
	is then given by the set $D^{(0)} = \{ x_L^{(0)}, \ldots,  x_R^{(0)} \}$
	with $x_L^{(0)} = \lfloor \min(w_1, \eta_1) \rfloor$
	and $x_R^{(0)} = \lceil \max(w_k, \eta_z) \rceil$.
	\newline
	We extend the support until 
	the relative change of the dual function 
	becomes smaller than the threshold $\dualFunctionThreshold$
	\begin{equation}
	\label{eq:dualFunctionRelativeChange}
	\left|\frac{\Psi(\lambda^{(\ell-1)}) -
	 \Psi(\lambda^{(\ell)})}{\Psi(\lambda^{(\ell)})}\right| < \dualFunctionThreshold,
	\end{equation}
	where we choose $\dualFunctionThreshold = 10^{-4}$
	for all case studies.
	If the inequality is not satisfied,
	we extend the support 
	by adding new states in each iteration
		\begin{equation}
			\hspace{-2.5ex} \left( x_L^{(\ell+1)}, x_R^{(\ell+1)} \right)\! =\! 
			\left( \max (0, x_L^{(\ell)} - 1), x_R^{(\ell)}+1 \right)
		\end{equation}
		The final results $\tilde{\lambda}$ and $\hat{D}$ of the iteration yields the distribution 
		$\tilde{q}(x)$ that approximates the marginal distribution of interest.


\subsection{Numerical Approach for the Two-dimensional Maximum Entropy Problem}
\label{appendix:maxEntProblem2D}
	In the case of two-dimensional distributions, 
	the maximum entropy problem is modified as follows. 
	We consider a sequence of non-central moments
	$\Ex{X_\centerdot^r X_\circ^l} = \mu_{r,l}$, $0 \leq r+l \leq M$,
	and the set $\distrFamily^2$ of all 
	two-dimensional discrete distributions 
	that satisfy the following constraints
	\begin{equation}
	\label{eq:momentconstr2d}
		\sum\limits_{x,y} x^r y^l g(x,y)   = \mu_{r,l}, \quad 0 \leq r+l \leq M.
	\end{equation}
	Here
	$X_\centerdot$ and $X_\circ$ correspond to the populations
	of two different species, i.e.\ to two distinct 
	elements of the random vector $\mcProcess(t)=(X_1(t),\ldots,X_n(t))$
	at some fixed time instant $t$.
	Similarly to the optimization problem~(3), 
	we seek 
	the 
	distribution $q \in \distrFamily^2$ that maximizes the entropy $H(g)$
	\begin{equation}
	\label{eq:maxShannonProblem2d}
	 \begin{array}{rcl}
		\maxent &= &\arg \max\limits_{g\in \distrFamily^2} H(g)\\[1ex]
				&=& \arg \max\limits_{g\in \distrFamily^2} 
				\left( -\sum\limits_{x,y} g(x,y) \ln{g(x,y)}\right)
	\end{array}
	\end{equation}
	We then proceed similarly to the one-dimensional case.
	 The general form of the solution 
	for the maximum entropy problem is given by
		\begin{equation}
		\label{eq:optSolution2D}
		\begin{array}{rcl}
		q(x,y) &= &
		\exp ( -1 -\sum\limits_{0 \leq r+l \leq M} \lambda_{r,l} x^r y^l ) \\[1ex]
		& =&\frac{1}{Z}
	\exp (-\sum\limits_{1 \leq r+l \leq M} \lambda_{r,l} x^r y^l ) ,
	\end{array}
	\end{equation}	
	where the normalization constant $Z$ is calculated as
	\begin{equation}
	\label{eq:normalizationConst2d}
	\textstyle Z = e^{1+\lambda_{0,0} }
     = \sum_{x,y} \exp 
     ( 
	      	-\sum\limits_{1 \leq r+l \leq M} \lambda_{r,l} x^r y^l
     ).
	\end{equation}
	

%
	We solve the optimization problem numerically similarly to the
	one-dimensional case.
	The vector
	$\lambda^{(\ell)} = 
	\left(
		\lambda_{0,1}, \lambda_{1,0}, \ldots, \lambda_{0,M}, \lambda_{M,0}
	\right)$
	is an approximation of the vector $\lambda$ in Eq.~\eqref{eq:optSolution2D}. 
	The elements of the gradient vector are computed as
	$ \nicefrac{\partial \Psi}{\partial \lambda_{r,l} } \approx
 \mu_{r,l} - (\nicefrac{1}{Z}) \widetilde{\mu}_{r,l}, $
	where  $\widetilde{\mu}_{r,l}$ is approximated by
	\begin{equation}
	\label{eq:momentApprox2D}
		\textstyle\widetilde{\mu}_{r,l} = 
		\sum_{x,y} x^r y^l \exp 
		 ( 
				-\sum\limits_{1 \leq r+l \leq M} \lambda_{r,l} x^r y^l.
		 ),
	\end{equation}
	Here $r,l \in\{ 0, \ldots, 2M\}$
	and the sum is taken over all $(x,y) \in \bbbn^2_0$.
	Finally, the elements of the Hessian matrix are computed as
	$$\textstyle H_{r+u,l+v} = 
	\frac{\partial^2 \Psi}{\partial \lambda_{r,l} \partial \lambda_{u,v}} 
	\approx
	\frac{Z \cdot \widetilde{\mu}_{r+u,l+v} - \widetilde{\mu}_{r,l} \widetilde{\mu}_{u,v}}{Z^2},
	$$
	where $0 \leq r+l \leq M, 0 \leq u+v \leq M$.
	Following the same procedure as in 
	Section~4.1, 
	the vector
	$\lambda^* = \left(
					\lambda^*_{0,1}, \lambda^*_{1,0}, \ldots, \lambda^*_{0,M}, \lambda^*_{M,0}
	\right)$ is found.
	The dimensionality of the optimization problem is 
	$0.5 \left(M^2 + 3M\right)$,
	and $\lambda^*_{0,0}$ can be calculated from~\eqref{eq:normalizationConst2d}
	as $\lambda^*_{0,0} = \ln{Z} - 1$.
	In comparison to the one-dimensional case,
	the range of the values of $\widetilde{\mu}_{r,l}$ becomes 
	wider due to the larger dimensionality, 
	so that the conditional number of the matrix
	$ \left( H + \gamma^{(\ell)} \cdot \mathrm{diag}(H) \right) $
	is even higher and the iteration might fail.
	\newline
	To approximate the moment values in~\eqref{eq:momentApprox2D}
	we truncate the infinite support 
	and consider the subset  $D_{xy}^{*} = D_{x}^{*} \times D_{y}^{*} $ instead. 
	Again, we choose $D_{xy}^{*} $ such that 
	the relative change of the dual function
	~\eqref{eq:dualFunctionRelativeChange}
	becomes smaller than the threshold $\dualFunctionThreshold$.
	The approximation $\tilde{q}(x,y)$ 
	of the marginal distribution 
	$p_{\centerdot,\circ}(x,y,t) = 
	P \left( X_\centerdot(t) = x, X_\circ(t) = y \right)$
	is then defined by the result $\lambda^*$ of the iteration procedure
	such that
	$p_{\centerdot,\circ}(x,y,t) \approx 
	\tilde{q}(x,y) \mbox{ if } (x,y) \in D_{xy}^{*}$
	and
	$p_{\centerdot,\circ}(x,y,t) \approx  0 \mbox { if } (x,y) \notin D_{xy}^{*}.$

\subsection{Reconstruction of Distributions from Approximated Moments}
\label{appendix:reconstructionMethods}
	In the following we discuss the details 
	of the reconstruction of marginal probability distributions  
	based on solving the moment problem using the maximum entropy approach.
	We consider the three possibilities introduced in Section~5, 
	\wsMCM{}, \jMCM{} and \MM{}.
	We illustrate the details of all three approaches
	with examples.
	\begin{exmp}
 	\label{ex:reconstructionMethodsDetails}
		We consider the gene expression model 
		(cf. Example 1) 
		where we reconstruct the marginal distribution of protein molecules
		$P \left( X_P(t) = x \right) = p_{X_P} (x,t) $.
		The  moments $\mu_k=\Ex{X_P^k}$ and the corresponding conditional
		moments are obtained   using the
		MCM and MM equations, for $k=0,\ldots,M+1$.
		In the case of \jMCM{} and \MM{}
		we use the first $M$ moments' values   as constraints in
		Eq.~(2) 
		and solve the maximum entropy optimization problem
		in Eq.~(3). 
		In both cases, the solution is given by a pair 
		$\left( \lambda^*, D^* \right)$
		of the parameter vector $\lambda^*$ and the truncated support $D^*$.
		The corresponding reconstructed distribution is defined as
		\begin{equation*}
		\begin{aligned}
			\tilde{q}(x) &=
				\exp \left(\!-\!1 \!-\! \sum\limits_{k=0}^{M} \lambda^*_{k} x^k \right)
				,  x \in D^*. \\
		\end{aligned}
		\end{equation*}
		In order to apply the \wsMCM{}, 
		we reconstruct the conditional distribution from the sequences 
		$ \mu_{P_{\mbox{\scriptsize off}},k} $ 
		and 
		$ \mu_{P_{\mbox{\scriptsize on}},k} $ 
		that approximate the conditional moments
		$\Ex{X_P^k | D_{\mbox{\scriptsize off}} =1}$
		and
		$\Ex{X_P^k | D_{\mbox{\scriptsize on}} =1}$
		.
		Here, $X_P$ corresponds to the number of proteins and
		the condition  $D_{\mbox{\scriptsize off}} =1$ 
		($D_{\mbox{\scriptsize on}} =1$) refers to the 
		state of the gene.
		These sequences of moments are obtained using the MCM approach together
		with the approximation of the mode probabilities 
		$p_{\mbox{\scriptsize off}}$ and $p_{\mbox{\scriptsize on}}$
		(cf.\ Example~\ref{ex:geneExpression_MCM}).
		We solve the maximum entropy  problem for each moment sequence
		and the reconstruction of marginal unconditional distribution 
		is given by 
		\begin{equation*}
				\tilde{q}_{wsMCM}(x) = 
				\begin{cases}
					p_{\mbox{\scriptsize off}} \tilde{q}_{\mbox{\scriptsize off}}(x), & 
						x \in D^*_{P_{\mbox{\scriptsize off}}} \setminus D^*_{P_{\mbox{\scriptsize on}}}\\
					p_{\mbox{\scriptsize on}} \tilde{q}_{\mbox{\scriptsize on}}(x), &
						x \in D^*_{P_{\mbox{\scriptsize on}}} \setminus D^*_{P_{\mbox{\scriptsize off}}}\\
					\begin{aligned}[c]
					p_{\mbox{\scriptsize off}} & \tilde{q}_{\mbox{\scriptsize off}}(x)  \\
					 + & p_{\mbox{\scriptsize on}} \tilde{q}_{\mbox{\scriptsize on}}(x)
					\end{aligned},
					 & x \in D^*_{P_{\mbox{\scriptsize off}}} \cap D^*_{P_{\mbox{\scriptsize on}}},
				\end{cases}
			\end{equation*}
		where 
		$\tilde{q}_{\mbox{\scriptsize off}}(x)$ and $\tilde{q}_{\mbox{\scriptsize on}}(x) $
		are the reconstructions of the conditional distributions. 
	\end{exmp}
	 
	To reconstruct two-dimensional marginal distributions 
	we numerically  solve the
	two-dimensional maximum entropy problem 
	as described in   Section~\ref{appendix:maxEntProblem2D}.
	We illustrate how   two-dimensional distributions are  reconstructed
	through
	the following example where we apply the \wsMCM{} approach.
	
	\stepcounter{figure}
	\stepcounter{figure}
	\stepcounter{figure}
	\stepcounter{figure}
	\stepcounter{figure}
	\begin{figure*}[th!]
		\centering
			\subfloat
			{
			\label{subfig:gexp_P_m3}
			\includegraphics[width=0.33\textwidth]
			{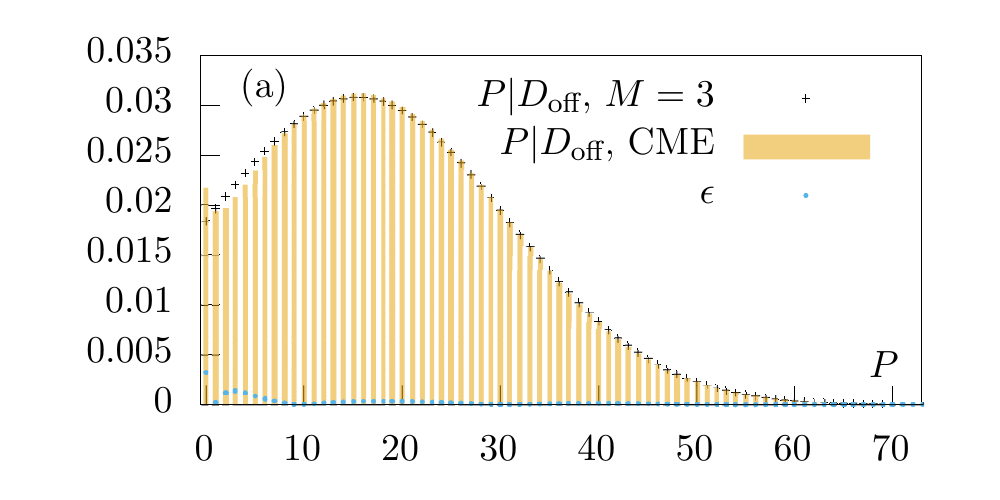}
			}
			\subfloat 
			{
			\label{subfig:gexp_P_m5}
			\includegraphics[width=0.33\textwidth]
			{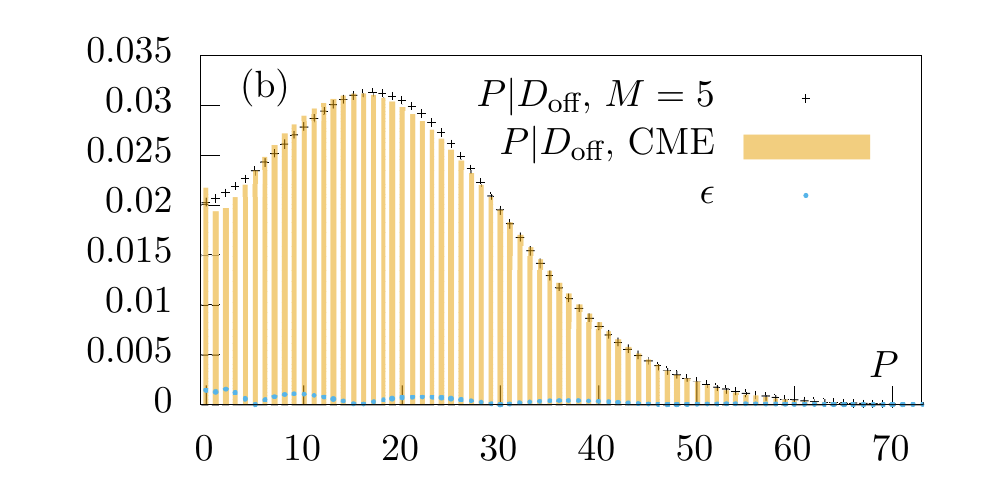}
			}
			\subfloat 
			{
			\label{subfig:gexp_P_m7}
			\includegraphics[width=0.33\textwidth]
			{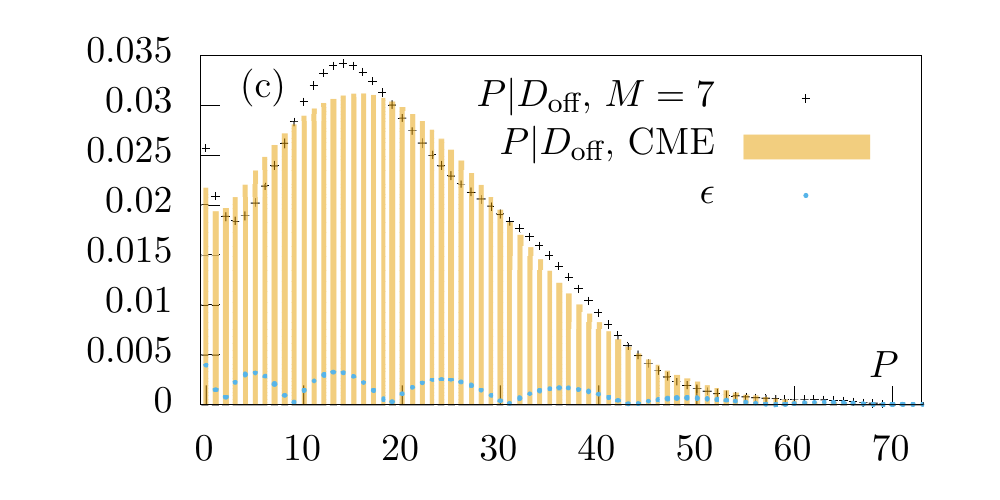}
			}
	 		\caption{\label{fig:P_Doff_comparison_different_M}
	 		Approximation of the conditional distribution of protein
	 		$P|D_{\mbox{\scriptsize off}}$.
	 		The number of moments used for reconstruction is
	 		$M=3$ (a),
	 		$M=5$ (b)
	 		and
	 		$M=7$ (c).}
\end{figure*}
	
	\begin{exmp}
	\label{ex:reconstructionMethodsDetails2D}
		We consider the exclusive switch system described 
		in Section~5. 
		The goal here is to reconstruct the 
		two-dimensional marginal distribution \hfill \\
		$P( X_{P_1} = x, X_{P_2}(t) = y )$
		of proteins $P_1$ and $P_2$.
		We first
		approximate
	  the mode probabilities  
		$p_1 = P \left({\DNA} = 1 \right)$,
		$p_2 = P \left({\DNA.P1} = 1 \right) $
		and
		$p_3 = P \left({\DNA.P2} = 1 \right)$
		(cf. Eq.~\ref{eq:smallcme}).
				In addition,  the conditional moments 
		\begin{equation*}
		\begin{aligned}
			\mu_{1;r,l} &= \Ex{X_{P_1}^r X_{P_2}^l | \DNA =1}  \\
			\mu_{2;r,l} &= \Ex{X_{P_1}^r X_{P_2}^l | \DNA.P1 =1}	\\
			\mu_{3;r,l} &= \Ex{X_{P_1}^r X_{P_2}^l | \DNA.P2 =1}
		\end{aligned}
		\end{equation*} are approximated 
		for $0 \leq r+l \leq M+1$,
		where $\DNA =1$ refers to the case where the promoter is free
		and $\DNA.P1 =1$ ($\DNA.P2 =1$) to the case where a molecule of type $P_1$
		(type $P_2$)
		is bound to the promoter.
		The constraints~\eqref{eq:momentconstr2d} for the 
		maximum entropy problem are given by the elements
		of these three sequences for $0 \leq r+l \leq M$ and
		the corresponding solutions of the optimization problem
		are given by the pairs 
		$\left( \lambda^*_i, D^*_i \right)$, $i=\{1,2,3\}$.
		Then the reconstructed  distribution
		is given by
		\begin{equation*}
		\begin{aligned}[c]
			& \tilde{q}_{wsMCM}(x,y) = \\
			& \begin{cases}
				p_1 \tilde{q}_1(x,y), & 
					(x,y) \in D^*_1 \setminus \left(D^*_2 \cup D^*_3\right) \\
				p_2 \tilde{q}_2(x,y), & 
					(x,y) \in D^*_2 \setminus \left(D^*_1 \cup D^*_2\right) \\
				p_3 \tilde{q}_3(x,y), & 
					(x,y) \in D^*_3 \setminus \left(D^*_1 \cup D^*_2\right) \\
				\sum\nolimits_{i=1}^2 p_i \tilde{q}_i(x,y), & 
					(x,y) \in \left(D^*_1 \cap D^*_2 \right) \setminus D^*_3,\\
				\sum\nolimits_{i=\{1,3\}} p_i \tilde{q}_i(x,y), & 
					(x,y) \in \left(D^*_1 \cap D^*_3 \right) \setminus D^*_2,\\
				\sum\nolimits_{i=2}^3 p_i \tilde{q}_i(x,y), & 
					(x,y) \in \left(D^*_2 \cap D^*_3 \right) \setminus D^*_1,\\
				\sum\nolimits_{i=1}^3 p_i \tilde{q}_i(x,y), &
					(x,y) \in D^*_1 \cap D^*_2 \cap D^*_3,
			\end{cases}
		\end{aligned}
		\end{equation*}
		where $\tilde{q}_i(x,y) = \exp (-1 - \sum_{1 \leq r+l \leq M} 
						\lambda^*_{r,l} x^r y^l )$.
	\end{exmp}

\subsection{Case Studies}
\label{appendix:case_studies_more_details}
Here we present   detailed results  of the reconstruction
of the marginal distributions that were discussed in Sect.~5. 

	\paragraph*{Gene Expression Model.}
	We show the approximation error
	$\maxDevPercErr$ 
	for the reconstruction
	of both conditional and unconditional distributions
	for mRNA and protein
	in Table~\ref{tab:gexp_1D_RP_together}
	where we use the first parameter set.
	Here, the first two columns refer to the
	approximation error of the conditional distributions for protein (mRNA)
	denoted by 
	$P|D_{\mbox{\scriptsize off}}$ ($R|D_{\mbox{\scriptsize off}}$) 
	and
	$P|D_{\mbox{\scriptsize on}}$ ($R|D_{\mbox{\scriptsize on}}$).
	The last three columns refer to the reconstructions
	of the marginal distribution 
	obtained using \wsMCM{}, \jMCM{} and
	\MM{}, 	respectively.

	
	\begin{table}[b]
			\centering
			\caption{
			Approximation errors of mRNA and protein 
			distributions reconstruction
			for gene expression system
			(first parameter set).
			\label{tab:gexp_1D_RP_together}
			}
			\nprounddigits{1}
\begin{tabular}{Hcccccc}
		\hline
					metric
					& M
					& $P|D_{\mbox{\scriptsize off}}$
					& $P|D_{\mbox{\scriptsize on}}$
					& $P_{wsMCM}$
					& $P_{jMCM}$
					& $P_{MM}$
					\\ \hline 
	\multirow{3}{*}{$\maxDevPercErr$} 
	& 3 & \numprint{15.0559} & \numprint{26.7230} & \numprint{6.5425} & \numprint{10.5432} & \numprint{10.2246} \\
	& 5 & \numprint{7.9403} & \numprint{1.7476} & \numprint{4.2689} & \numprint{9.9950} & \numprint{8.1045} \\
	& 7 & \numprint{18.4318} & \numprint{2.7662} & \numprint{12.5378} & \numprint{7.9192} & \numprint{3.2601} \\
	\hline
				metric
				& M
				& $R|D_{\mbox{\scriptsize off}}$ 
				& $R|D_{\mbox{\scriptsize on}}$ 
				& $R_{wsMCM}$ 
				& $R_{jMCM}$ 
				& $R_{MM}$
				\\ \hline 
	\multirow{3}{*}{$\maxDevPercErr$} 
	& 3 & 1.0 & 19.9 & 1.2 & 3.9 & 3.8 \\
	& 5 & 0.5 & 4.6 & 1.1 & 0.2 & 0.4 \\
	& 7 & 1.3 & 0.6 & 1.3 & 0.3 & 0.6 \\ \hline

\end{tabular}
\npnoround

	\end{table}
	
	We observe that the reconstruction is most accurate 
	for the distribution of mRNA when the
 	\jMCM{} method is applied with $M=7$.
 	The distribution of protein molecules is 
 	reconstructed	most accurately 
 	when \MM{} is applied with $M=7$.
	Please note that the large approximation errors of 
	conditional distribution reconstructions may still provide 
	an accurate reconstruction for the unconditional distribution
	because of the computation of a weighted sum that can average out
	individual deviations from the true probability value.
	For example, the reconstruction of the marginal distribution of 
	proteins and mRNA with \wsMCM{} 
	gives the smallest error when we use $M=3$ moments
	(6.5\% for proteins, 1.3\% for mRNA),
	but the approximation errors of the corresponding conditional
	distributions are much larger.
	
	The sensitivity of the optimization procedure can also influence the 
	final result. 
	The reconstruction that uses fewer degrees of freedom 
	can provide an accurate solution since the distribution
	of the simple shape is able to explain the main behavior.
	At the same time, adding more moments 
	into the consideration allows one to capture more details,
	but it may change the reconstruction drastically
	due to the sensitivity, and the corresponding approximation error can become larger.
	To the best of our knowledge, there exist no criteria
	that provide the number of moments that have to be
	considered such that adding more information does not 
	greatly
	change
	the maximum entropy reconstruction.
	We show that in Fig.~\ref{fig:P_Doff_comparison_different_M}, 
	where we plot the reconstructions of the conditional
	distribution $P|D_{\mbox{\scriptsize off}}$
	and use $M \in \{3,5,7\}$.
	The reconstruction using $M=7$ moments has the largest approximation error,
	but it is able to capture the complex nature of the distribution
	by treating the point $P=0$ differently.

		\begin{figure*}[t]
	\centering
		\subfloat{
			\label{subfig:gexp_slow_1D_P}
			\includegraphics[width=0.5\textwidth] 
				{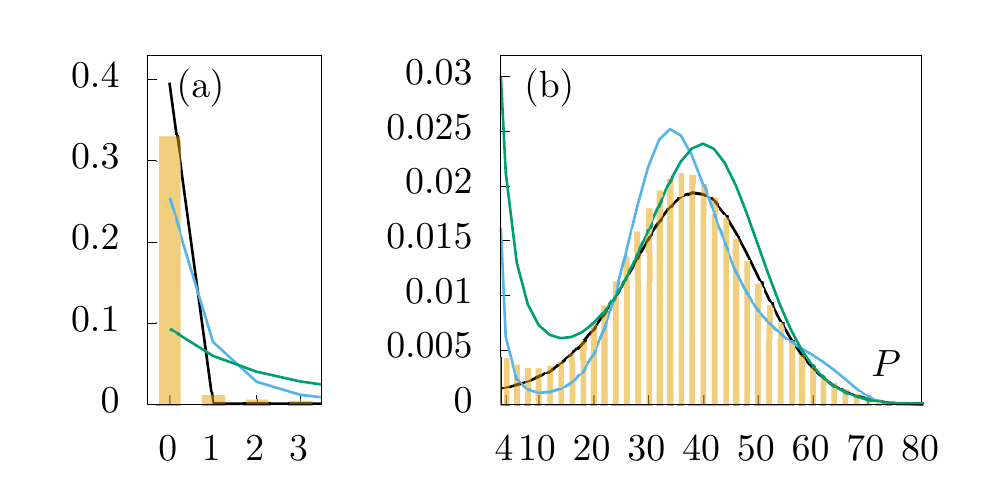}
			}
		\subfloat{
		\label{subfig:gexp_slow_1D_R}
			\includegraphics[width=0.5\textwidth] 
			{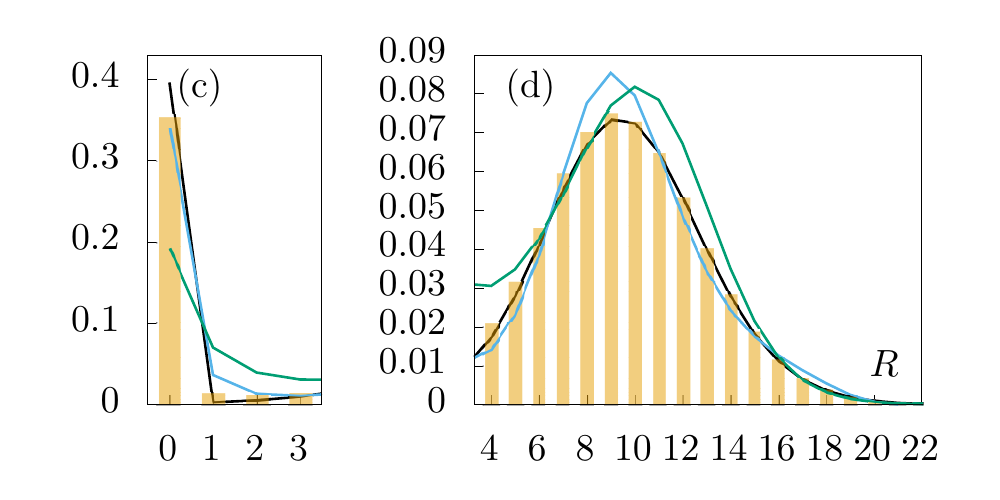}
		}
 		\caption{\label{fig:gexp_slow_R_P}
 		Gene expression (second parameter set): 
 		approximations of the marginal distribution of 
		protein (a), (b) 
 		and
		mRNA (c), (d) 
 		obtained using $M=5$ moments.
 		The reconstructions are plotted with lines
 		(black for \wsMCM{},
 		blue for \jMCM{} 
 		and green for \MM{})
 		and the CME solution is plotted with bars.
 		The plots (a) and (c) show in detail the region
 		with molecular counts $\{0,1,2,3\}$.
 		}
	\end{figure*}

 	We note that the reconstruction results are 
 	generally quite similar for the approaches
 	that are based on an approximation of the
	conditional and unconditional moments.
 	However, the MCM approach has the advantage 
 	that the distribution of species such as DNA 
 	is very accurate, since they are directly available 
 	and are not reconstructed
	from the moments.
	A moment-based approach such as MM needs a 
	large number of moments for an accurate reconstruction~\cite{bandyopadhyay2005maximum}.
	We also notice that the approximation of the conditional moments 
	in the MCM method
	is less	accurate than the approximation of the unconditional moments
	in the MM method (cf.\ Tables~\ref{tab:gexpr_mc_error} and~\ref{tab:gexpr_mcm_error})
	for this parameter set.
				\begin{table}[b!]
				\centering
				\caption
				{ Approximation errors of two-dimensional	distribution reconstruction
					for gene expression example (first parameter set).
					\label{tab:genexp_2D_reconstruction}	}
				{	

\nprounddigits{1}
\begin{tabular}{Hcccccc}
	\hline
				metric
				& M
				& $\tilde{q}_{D_{\mbox{\scriptsize off}}}$
				& $\tilde{q}_{D_{\mbox{\scriptsize on}}}$
				& $\tilde{q}_{wsMCM}$ 
				& $\tilde{q}_{jMCM}$ 
				& $\tilde{q}_{MM}$ 
				\\ \hline 
\multirow{3}{*}{$\maxDevPercErr$} 
& 3 & \np{59.1054} & \np{48.8522}	& \np{58.1141} & \np{62.1209} & \np{61.9539} \\ 
& 5 & \np{53.4683} & \np{36.1533}	& \np{51.7948} & \np{58.1554} & \np{58.0910}\\ 
& 7 & \np{28.0351} & \np{29.1008} & \np{28.1394} & \np{24.7305} & \np{29.135} \\ \hline
\end{tabular}
\npnoround }
				\end{table}
	Nevertheless, the 
		reconstruction based on the conditional moments is in some cases more
		accurate, which means that the error
		is mostly due to the maximum entropy procedure.
	
	An example of a two-dimensional distribution reconstruction
	is shown in Fig.~4. 
	Here we present in addition the approximation errors for all three reconstruction
	methods 
	in Tab.~\ref{tab:genexp_2D_reconstruction},
	both for conditional and marginal two-dimensional
	distributions of mRNA and protein.
	For the sake of readability we denote the reconstructed 
	distribution by $\tilde{q}$ in the following tables.
	For instance, the approximation of the joint marginal distribution of
	$R$ and $P$ under the condition $D_{\mbox{\scriptsize off}} = 1$
	is denoted by $\tilde{q}_{D_{\mbox{\scriptsize off}}}$.
	We observe that the approximation error decreases when we make
	use of more moments.

	\begin{table}[b]
			\centering
			\caption{
			Approximation errors of reconstructed mRNA and protein 
			distributions
			for the gene expression example 
			(second parameter set).
			\label{tab:gexp_slow_1D_RP_together}
			}
			\nprounddigits{1}
\begin{tabular}{Hcccccc}
		\hline
					metric
					& M
					& $P|D_{\mbox{\scriptsize off}}$
					& $P|D_{\mbox{\scriptsize on}}$
					& $P_{wsMCM}$
					& $P_{jMCM}$
					& $P_{MM}$
					\\ \hline 
	\multirow{3}{*}{$\maxDevPercErr$} 
	& 3 & 9.5 & 93.0 & 8.5 & 59.8 & 88.9 \\
	& 5 & 21.3 & 70.3 & 20.1 & 23.1 & 71.6 \\
	& 7 & 21.3 & 78.4 & 20.0 & $>100$ & 60.7 \\
	\hline
				metric
				& M
				& $R|D_{\mbox{\scriptsize off}}$ 
				& $R|D_{\mbox{\scriptsize on}}$ 
				& $R_{wsMCM}$ 
				& $R_{jMCM}$ 
				& $R_{MM}$
				\\
				\hline 
	\multirow{3}{*}{$\maxDevPercErr$} 
	& 3 & $>100$ & 10.7 & 85.9 & 25.1 & 71.5 \\
	& 5 & 12.4 & 2.5 & 12.1 & $>100$ & 45.6 \\
	& 7 & 12.4 & 1.3 & 12.2 & 46.1 & 33.7 \\
	\hline
\end{tabular}
\npnoround

	\end{table}
	
	The computation time for the reconstruction of the
	one-dimensional distribution 
	for our MATLAB implementation
	(on a machine with the quad-core processor, 1.60GHz
	and 12 GB of RAM memory)
	is up to $0.3$ seconds
	whereas the approximation of the two-dimensional distribution 
	takes up to $15$ seconds.
	The running time mainly depends on the support approximation
	procedure.

	Next we consider the second parameter set.
	The approximation errors of the one-dimensional distributions
	are given in Table~\ref{tab:gexp_slow_1D_RP_together}.
	The results of the reconstruction for this parameter set
	are worse than for the first one
	due to the more complex shape of the distribution.
	The \wsMCM{} provides the best result
	in all the cases (except for the reconstruction
	of the mRNA distribution when $M=3$ moments are used).
	In Fig.~\ref{fig:gexp_slow_R_P}
	we show the reconstructions 
	both for protein and mRNA distribution
	obtained using $M=5$ moments.
	It can be seen that the maximum relative error $\maxDevPercErr$
	does not optimally describe the difference
	between the distribution shapes.
	For instance, a visual comparison
	of the reconstructed mRNA distributions
	reveals that the results obtained with \jMCM{} 
	describes the	shape better than the \MM{} based reconstruction
	though the relative error of \jMCM{} ($>100\%$)
	is larger than that of \MM{} ($45.6\%$).
	
	The approximation errors of the two-dimensional distributions
	are given in Table~\ref{tab:gexp_slow_2D}.
	It can be seen that the  results
	are worse than those of the first parameter set
	and taking more moments into consideration
	does not give better results.
	Thus,  the   entropy maximization 
	may not the best choice for the reconstruction
	of bi-modal distributions  where the values of the 
	peaks are of different orders of magnitude.

	\begin{table}[b]
			\centering
			\caption{
			Approximation errors of two-dimensional 
			reconstruction for gene expression example
			(second parameter set).
			\label{tab:gexp_slow_2D}
			}
			\nprounddigits{1}
\begin{tabular}{Hcccccc}
	\hline
				metric
				& M
				& $\tilde{q}_{D_{\mbox{\scriptsize off}}}$
				& $\tilde{q}_{D_{\mbox{\scriptsize on}}}$
				& $\tilde{q}_{wsMCM}$ 
				& $\tilde{q}_{jMCM}$ 
				& $\tilde{q}_{MM}$ 
				\\ \hline 
\multirow{3}{*}{$\maxDevPercErr$} 
& 3 & \np{82.6079} & \np{98.0822} & \np{82.6358} & \np{74.9253} & \np{95.2128} \\
& 5 & \np{83.6442} & \np{83.4017} & \np{83.5908} & \np{76.0199} & \np{88.0619} \\ 
& 7 & \np{92.4953} & \np{91.8545} & \np{92.4919} & \np{86.7912} & \np{89.9848} \\ \hline
\end{tabular}
\npnoround
	\end{table}

	
\begin{figure*}
	\centering
	\subfloat{
		\label{fig:exsw_m6_P1_DNA}
		\includegraphics[width=0.33\textwidth]
		{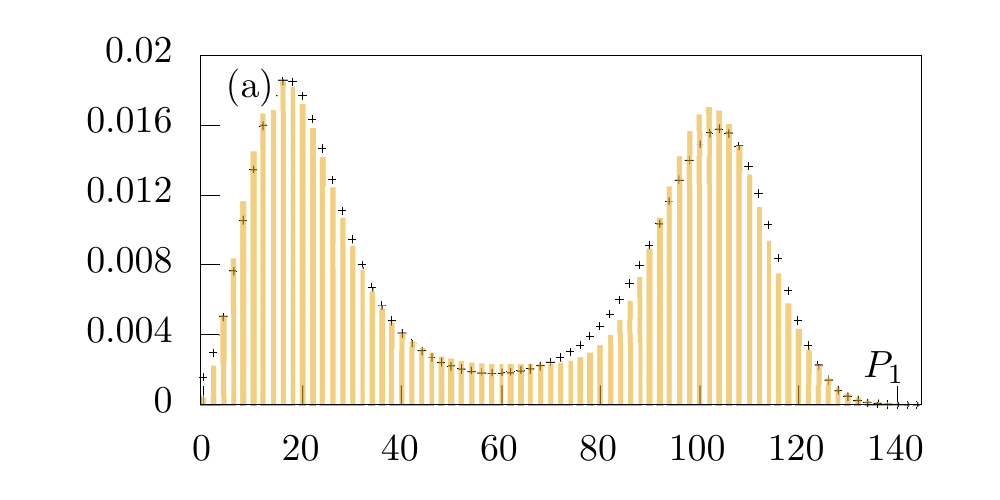}}
	\subfloat{
		\label{fig:exsw_m6_P1_DNA_P1_and_P1_DNA_P2}
		\includegraphics[width=0.33\textwidth]
		{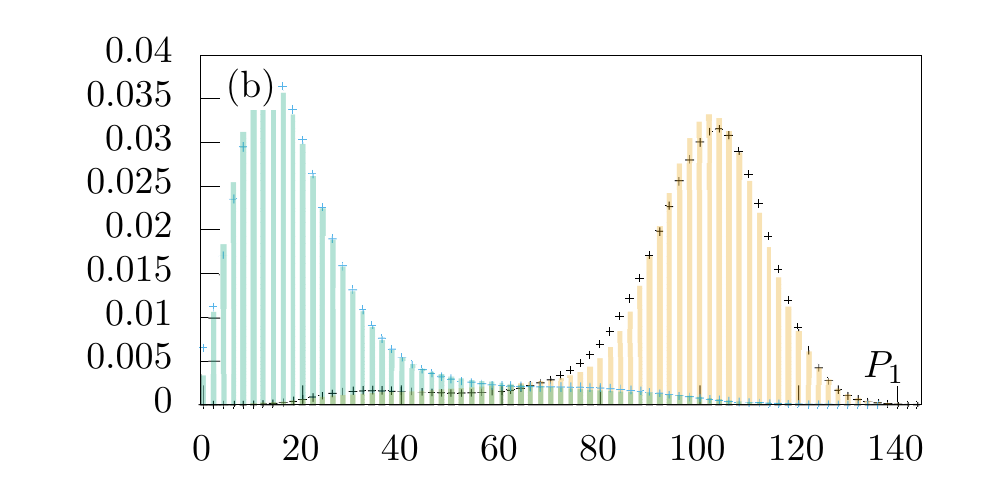}}
	\subfloat{
		\label{fig:exsw_m6_P1_wsMCM}
		\includegraphics[width=0.33\textwidth]
		{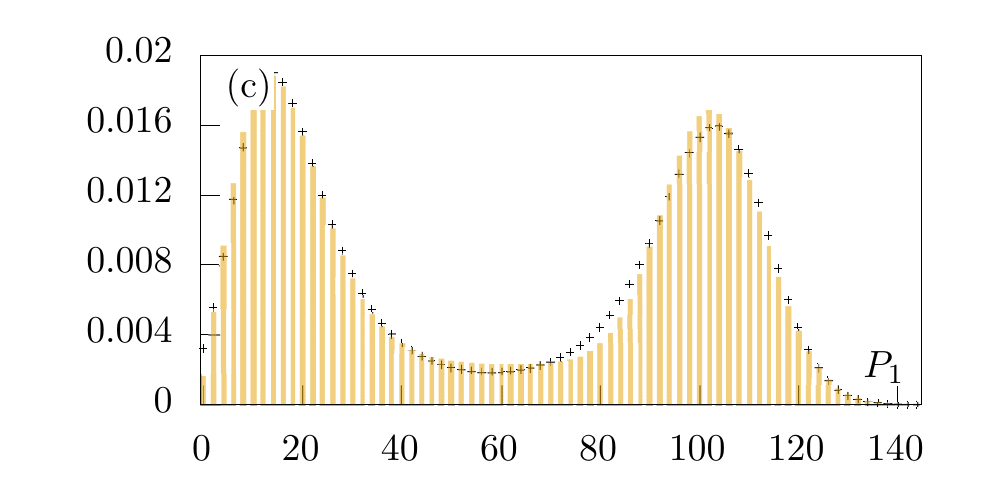}}
	\caption{\label{fig:exsw_m6_P1_reconstructions_plus_wsMCM}
	Exclusive switch: approximations of 
	the conditional distributions of protein $P_1$ where
	$DNA=1$ (a),
	$DNA.P_1=1$, $DNA.P_2=1$ (b)
	and
	the reconstruction of the marginal distribution (c).
	The solution of the CME is plotted with yellow bars 
	and
	the reconstructions are plotted with black crosses 
	(green bars and blue crosses are used for the conditional distribution (b) where $DNA.P_2=1$).
	The reconstruction of the marginal distribution (c) 
	is obtained using \jMCM{} with $M=5$.}
\end{figure*}

	\paragraph*{Exclusive Switch Model.}
	Next we address the  accuracy of the reconstruction
	of conditional and marginal distributions
	of the exclusive switch model introduced
	in Sect.~5. 
	In Table~\ref{tab:exswitch_tomacs_P1_P2_together}
	the approximation errors are listed for the conditional distributions of the proteins
	where we condition on the three possible states of the promoter, i.e.,
	 $\DNA=1$, $\DNA.P_1=1$ or $\DNA.P_2=1$.
	
		\begin{table}[b]
			\centering
			\caption
			{
				Approximation errors 
				for the distribution of 
				proteins $P_1$ and $P_2$.
				\label{tab:exswitch_tomacs_P1_P2_together}
			}
			{
\nprounddigits{1}
\begin{tabular}{cccccccc}
	\cline{2-8}
				& M
				& $\tilde{q}_{\DNA}$
				& $\tilde{q}_{\DNA.P_1}$
				& $\tilde{q}_{\DNA.P_2}$
				& $\tilde{q}_{wsMCM}$
				& $\tilde{q}_{jMCM}$
				& $\tilde{q}_{MM}$
				\\ 
				\cline{2-8}
\multirow{3}{*}{$P_1$} 
 & 3 & $>$100 & 20.7 & $>$100 & $>$100 & $>$100 & $>$100 \\
& 5 & 10.7 & 7.1 & 82.7 & 84.3 & $>100$ & $>100$ \\
& 7 & 5.7 & 7.8 & 79.5 & 80.5 & 6.6 & $>$100 \\
\cline{2-8}
\multirow{3}{*}{$P_2$} 
  & 3 & $>$100 & $>$100 & 40.8 & 41.5 & $>$100 & $>$100 \\
  & 5 & 17.7 & $>$100 & 14.9 & 14.8 & 17.8 & $>$100 \\
  & 7 & 15.3 & 7.4 & 7.5 & 8.2 & 12.0 & 19.1 \\
\cline{2-8}
\end{tabular}
			}
		\end{table}
	
	We observe that the 
	approximation error $\maxDevPercErr$ is minimal
	for both proteins $P_1$ and $P_2$
	when the \wsMCM{} approach is applied 
	for all $M \in \{3,5,7\}$.
	Thus, for the exclusive switch system 
	it is advantageous  to approximate the marginal distributions 
	by first reconstructing the conditional
	distributions  and computing the weighted sum.
	In almost  all cases the error  decreases when 
	more information about the moments is used.
	Because of the complex bi-modal shape of the distributions,
	it is beneficial to consider higher-order moments.
	It is important to note also that the large value 
	of the error ($\maxDevPercErr > 100$)
	usually comes from the
	probabilities around the boundary points of the support ($x_L$ or $x_R$).
	In the remaining parts of the support $D^{*}$ the reconstruction is accurate.
	For example, in Fig.~\ref{fig:exsw_m6_P1_reconstructions_plus_wsMCM}
	we show the reconstructions
	of both 
	conditional 
	(left and middle plots)
	and 
	marginal (right plot)
	distributions of $P_1$.
	Here, the \jMCM{} was used with $M=5$ to reconstruct 
	the marginal distribution.
	The visual comparison reveals 
	that the approximation nicely describes the bi-modal
	shape
	although the maximum relative error is large ($\maxDevPercErr>100$).

	We also consider the conditional and marginal
	two-dimensional distributions of proteins $P_1$ and $P_2$
	in Table~\ref{tab:exswitch_2d_reconstruction}.
	Again we condition on the state of the promoter
	region, e.g. $\tilde{q}_{\DNA.P_1}$ corresponds to the 
	joint distribution of proteins $P_1$ and $P_2$
	when $\DNA.P_1=1$.
	
	\begin{table}[b]
		\centering
		\caption{
		Two-dimensional conditional protein distributions (exclusive switch).
		\label{tab:exswitch_2d_reconstruction}
		}

\begin{tabular}{ccccccc}
		\hline
			M
			& $\tilde{q}_{\DNA}$
			& $\tilde{q}_{\DNA.P_1}$
			& $\tilde{q}_{\DNA.P_2}$
			& $\tilde{q}_{wsMCM}$
			& $\tilde{q}_{jMCM}$
			& $\tilde{q}_{MM}$
			\\ \hline 
		3 & 69.089 & $>100$ & 49.3173
			&  53.8828 	& $>100$ & $>100$ \\ 
		5 & 32.5258	& $>100$		& 47.1483
			&  45.5285			& 24.4008			& $>100$ \\ 
		7 & 19.7620 & $>100$   & 12.9727
			&  14.2448    & 28.2168    & 26.2527 \\ \hline
\end{tabular}

	\end{table}
	
\begin{figure*}
		\centering
		\subfloat{
			\label{fig:exsw_p1_p2_dna_5m}
				\includegraphics[width=0.33\textwidth]
				{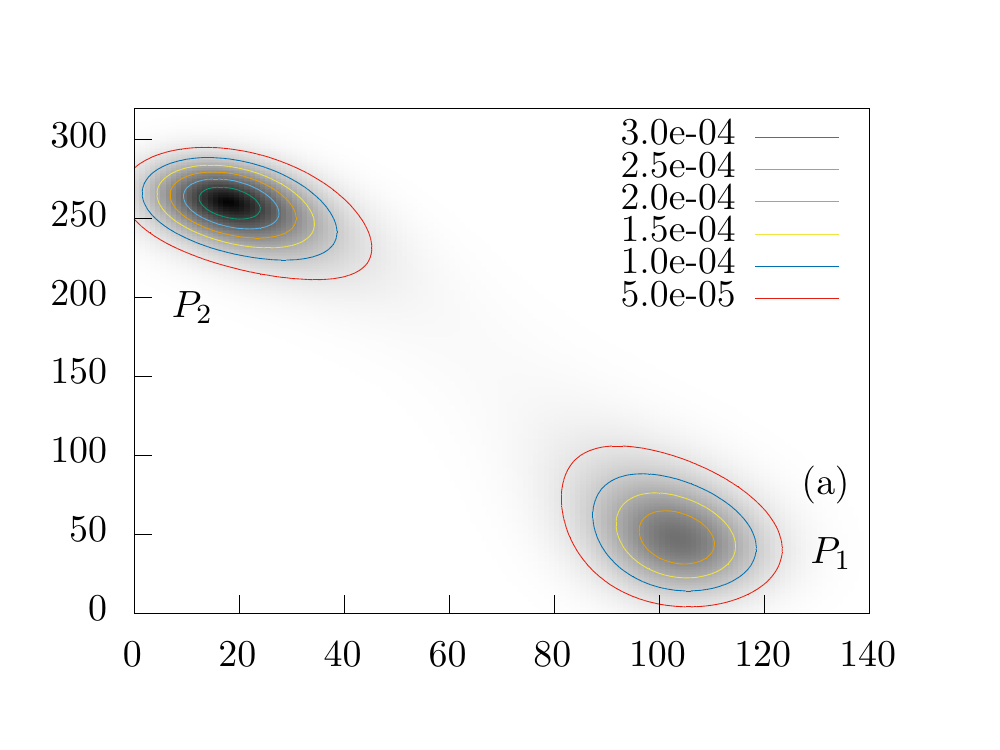}}
		\subfloat{
			\label{fig:exsw_p1_p2_dnaP1_5m}
				\includegraphics[width=0.33\textwidth]
				{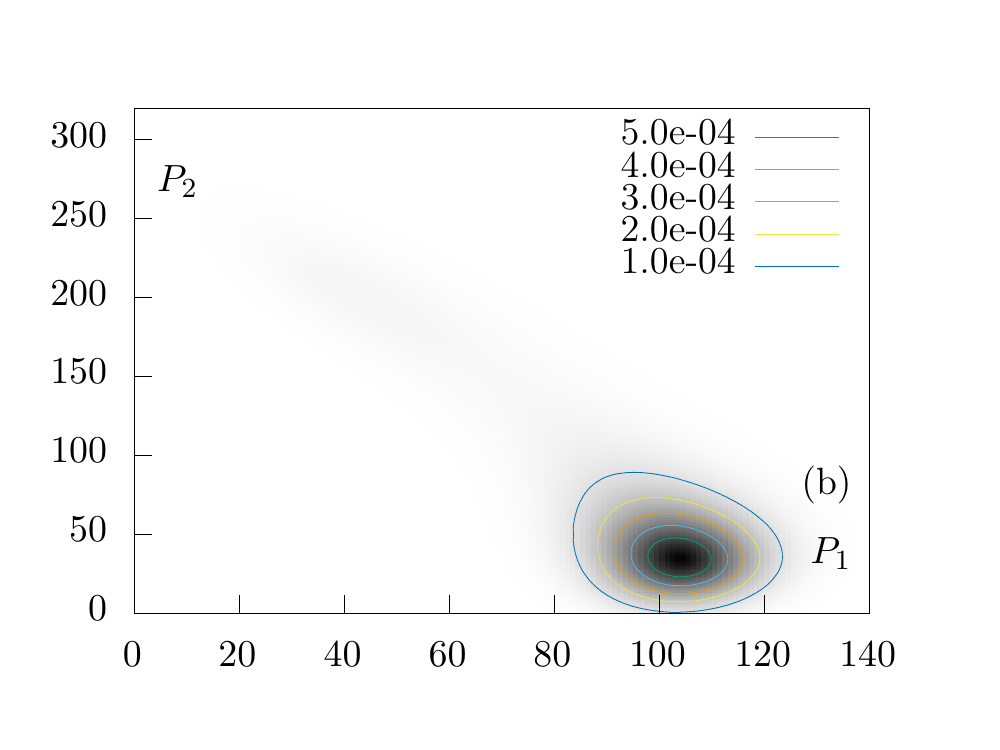}}
		\subfloat
			{\label{fig:exsw_p1_p2_dnaP2_5m}
				\includegraphics[width=0.33\textwidth]
				{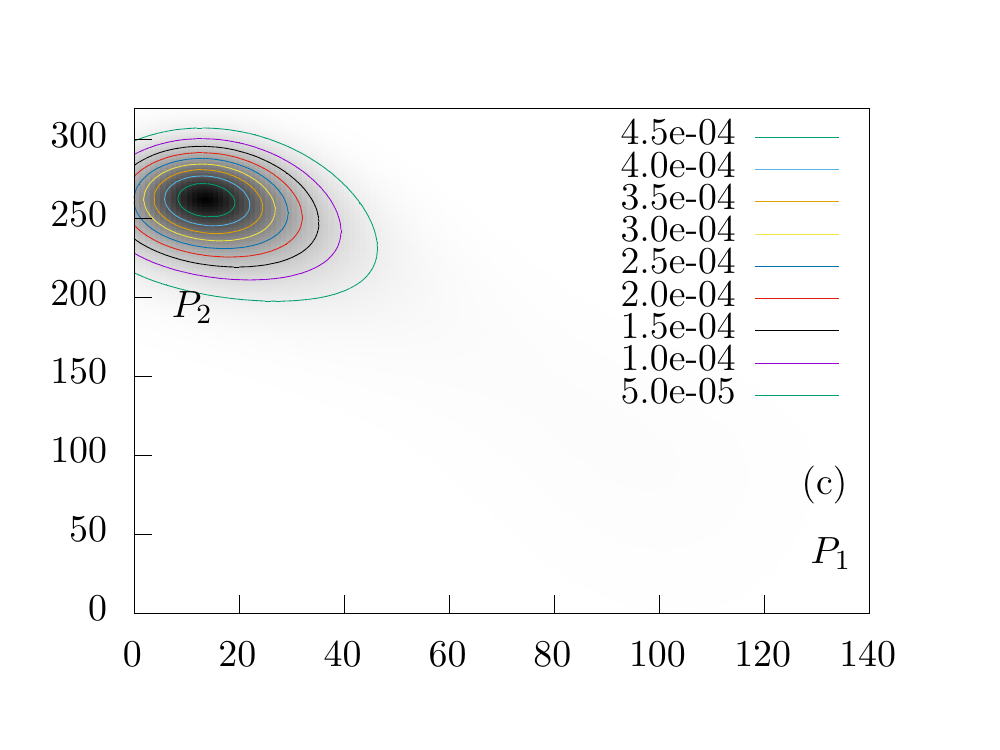}}
		\caption{\label{fig:exsw_2d_conditional_m5}
		Exclusive switch: the approximations of the conditional distributions
		of proteins $P_1$ and $P_2$
		where 
		$DNA=1$ (a),
		$DNA.P_1 = 1$ (b)
		and
		$DNA.P_2 = 1$ (c).
		The reconstructions are obtained using $M=5$ moments.}
	\centering
			\subfloat{\label{fig:exsw_p1_p2_unc_apprx}
				\includegraphics[width=0.33\textwidth]
				{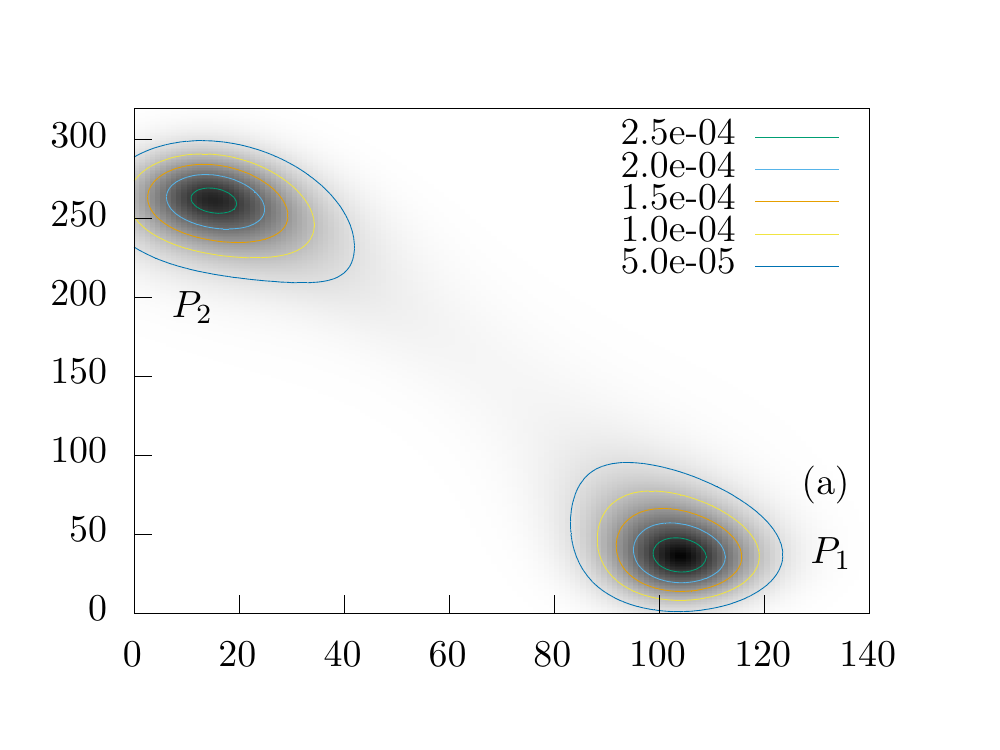}
				}
	\subfloat{\label{fig:exsw_p1_p2_unc_5m_distanceMap}
				\includegraphics[width=0.33\textwidth]
				{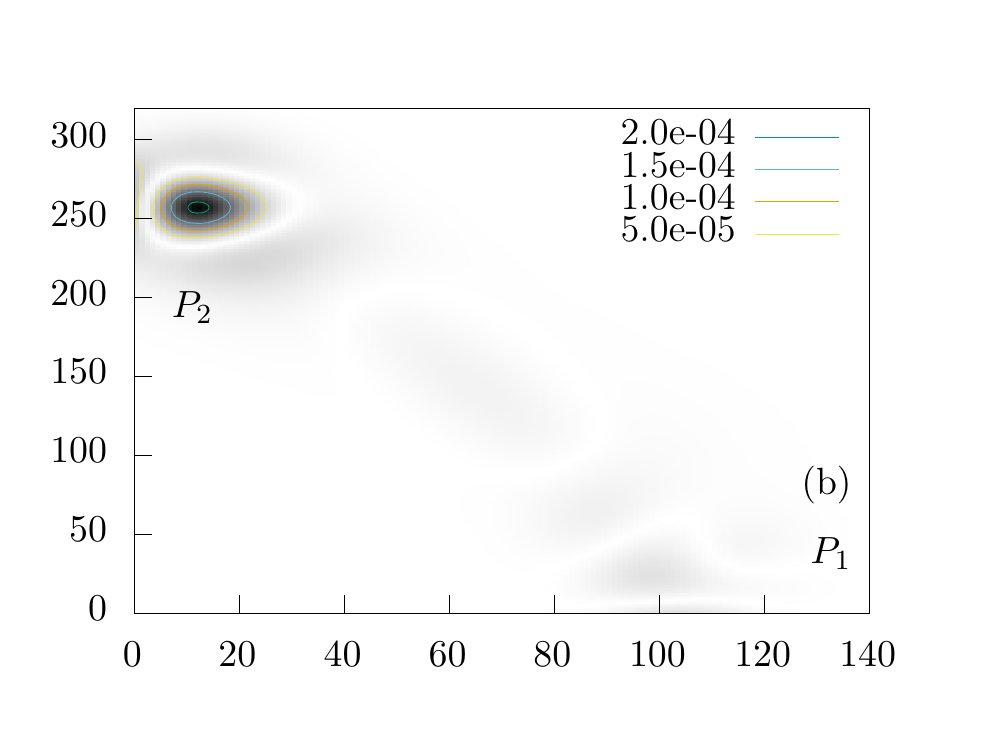}}
	\subfloat{\label{fig:exsw_5m_distanceMap_P2_slice}
				\includegraphics[width=0.33\textwidth]
				{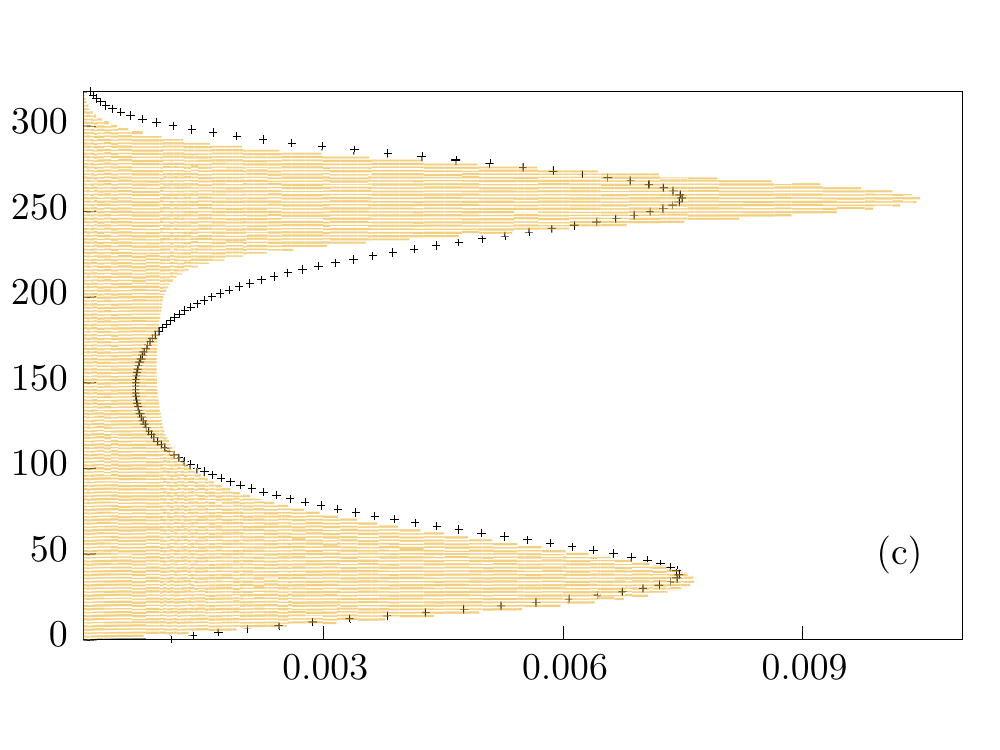}}
	\caption{\label{fig:exsw_2d_unc_m5}
	Exclusive switch:
	the reconstruction of the marginal distribution of proteins $P_1$ and $P_2$ (a)
	obtained using \wsMCM{} with $M=5$ moments
	and the corresponding approxmation error (b).
	The one-dimensional marginal distribution of $P_2$ (yellow bars)
	and the reconstruction (black dots) are shown in (c).}
	\end{figure*}

	The marginal distribution $P \left( X_{P_1}=x, X_{P_2}=y \right)$
	is best approximated when the \emph{weighted sum MCM} approach is applied
	and we   see that better reconstructions are achieved with higher
	order moments.
	Generally, the MCM approach gives more accurate results, i.e.,
	both \wsMCM{}
	and \jMCM{} perform better than \MM{}.
	We show the reconstructions of three conditional distributions
	in Fig.~\ref{fig:exsw_2d_conditional_m5} for the case when $M=5$,
	where the plots refer to the conditions
	(from left to right)
	$DNA=1$, $DNA.P_1=1$ and $DNA.P_2 = 1$.
	The reconstruction of the marginal distribution
	obtained using \wsMCM{}
	is shown 
	together with the approximation error
	in Fig.~\ref{fig:exsw_2d_conditional_m5}
	(left and right plot).
	We observe that the approximation error is
	large in this case.
	In Fig.~\ref{fig:exsw_2d_conditional_m5} 
	we also plot the marginal distribution of $P_2$ 
	where the mismatch for 
	the first peak of the distribution 
	can be explicitly seen.
	The reconstruction process for the exclusive switch model 
	takes more time 
	than for  gene expression model
	because of a much larger support.
	The running time of the one-dimensional reconstruction   
	  is up to $3$ seconds
	and in the two-dimensional case is up to $5$ minutes.
	Again, here the bottleneck of the reconstruction procedure is the support approximation.\\
	Thus, the idea of decomposing the Markov process
	into two  parts, as done for the conditional moment equations, results in 
	fewer equations and a more accurate description of the process.
	The weighted sum of mode probabilities and reconstructed 
	 conditional distributions seems to be particularly 
 	 beneficial when systems exhibit complex behavior, such as in the
	  exclusive switch  model.
	  \newpage
	


 \end{appendix}
 \bibliographystyle{plainnat}
 \bibliography{hsb}

\end{document}